%% file: acces.tex
\documentclass[12pt]{amsart}%

\usepackage{psfrag}
\usepackage{mathrsfs}
\usepackage{amsmath,amsthm}
\font\de=cmssi12
\usepackage{amssymb}
\usepackage{amsfonts}
\usepackage[dvips]{graphicx}
\usepackage{bibentry}
\usepackage{enumerate}
\usepackage{euscript}
\usepackage{yfonts}

\begin{document}


\newcounter{theorem}[section]
\newtheorem{conj}{\sc Conjecture}
\newcounter{\theconj}[conj]
\newtheorem{defi}[theorem]{\sc Definition}
\newtheorem{lema}[theorem]{\sc Lemma}
\newtheorem{proposition}[theorem]{\sc Proposition}
\newtheorem{corollary}[theorem]{\sc Corrollary}
\newtheorem{theorem}[theorem]{\sc Theorem}
\newtheorem{obs}[theorem]{\sc Remark}
\newtheorem{ques}[theorem]{\sc Question}
\newtheorem{question}{\sc Question}
\def\bp{\noindent{\it Proof. }}
\def\ep{\noindent{\hfill $\fbox{\,}$}\medskip\newline}
\renewcommand{\theequation}{\arabic{section}.\arabic{equation}}
\renewcommand{\thetheorem}{\arabic{section}.\arabic{theorem}}
\newcommand{\eps}{\varepsilon}
\newcommand{\disp}[1]{\displaystyle{\mathstrut#1}}
\newcommand{\fra}[2]{\displaystyle\frac{\mathstrut#1}{\mathstrut#2}}
\newcommand{\dif}{{\rm Diff}}
\newcommand{\phr}[1]{{\rm Par}^{#1}_m}
\newcommand{\Per}{{\mathcal P}er}
\newcommand{\Z}{\mathbb Z}
\newcommand{\R}{\mathbb R}
\newcommand{\N}{\mathbb N}
\newcommand{\T}{\mathbb T}
\newcommand{\s}{\sigma}
\newcommand{\W}{\mathcal W}
\newcommand{\F}{\mathcal F}
\newcommand{\G}{\mathcal G}
\newcommand{\He}{\mathcal H}
\newcommand{\go}[1]{\mathfrak {#1}}
\newcommand{\he}{\go{h}}
\newcommand{\length}{{\rm length}}
\newcommand{\diam}{{\rm diam}}
\newcommand{\vol}{\mathsf v}
\renewcommand{\j}{\jmath}
\def\to{\mathop{\rightarrow}}
\def\ord{\mathop{\rm ord}}
\def\diam{\mathop{\rm diam}}
\def\cc{\mathop{\rm cc}}
\title{Partial hyperbolicity and ergodicity in dimension three.}
\author{F.Rodriguez Hertz}
\author{M.Rodriguez Hertz}
\author{R. Ures}
\thanks{This work was partially supported by CONICYT-PDT 29/220 and CONICYT-PDT 54/18 grants}
\address{IMERL-Facultad de Ingenier\'\i a\\ Universidad de la
Rep\'ublica\\ CC 30 Montevideo, Uruguay} \email{frhertz@fing.edu.uy}
\email{jana@fing.edu.uy} \email{ures@fing.edu.uy}
\date{\today}
\begin{abstract} In \cite{rru} the authors proved the Pugh-Shub
conjecture for partially hyperbolic diffeomorphisms with
1-dimensional center, i.e. stable ergodic diffeomorphism are dense
among the partially hyperbolic ones. In this work we address the
issue of giving a more accurate description of this abundance of
ergodicity. In particular, we give the first examples of manifolds
in which {\em all} conservative partially hyperbolic diffeomorphisms
are ergodic.
\end{abstract} \maketitle
\input{intro2311.tex}
\input{prelim2311.tex}
\input{lami2311.tex}
\input{invtori2311}
\input{growth2311.tex}
\input{heisen2311.tex}
\input{nilergodic2311}


\end{document}

%% file: intro2311.tex
\begin{section}{Introduction}

A diffeomorphism $f:M\rightarrow M$ of a closed smooth manifold $M$
is partially hyperbolic if $TM$ splits into three invariant bundles
such that two of them, the strong bundles, are hyperbolic (one is
contracting and the other expanding) and the third, the center
bundle, has an intermediate behavior (see the next section for a
precise definition). An important tool in proving ergodicity of
these systems is the accessibility property i.e. $f$ is accessible
if any two points of $M$ can be joined by a curve that is a finite
union of arcs tangent to the strong bundles.

 In the last years, since the pioneer work of Grayson, Pugh and
Shub \cite{gps}, many advances have been made in the ergodic theory
of partially hyperbolic diffeomorphisms. In particular, we want to
mention the very recent works: of Burns and Wilkinson \cite{bw}
proving that (essential) accessibility plus a bunching condition
(trivially satisfied if center bundle is one dimensional) implies
ergodicity and of the authors \cite{rru} obtaining the Pugh-Shub
conjecture about density of stable ergodicity for conservative
partially hyperbolic diffeomorphisms with one dimensional center.
That is, ergodic diffeomorphisms contain an open and dense subset of
the conservative partially hyperbolic ones. See \cite{rru2} for a
recent survey on the subject.

After these result a question naturally arises:
 Can we describe this abundance of
ergodicity  more accurately? \par
More precisely:

\begin{ques}
Which manifolds support a non-ergodic partially hyperbolic
diffeomorphism?
\end{ques}

Somewhat surprisingly, studying this question we realize that, in
dimension 3, there are strong obstructions for a partially
hyperbolic diffeomorphism to fail to be ergodic.  We conjecture that
the answer to this question, in dimension 3, is that the only
(orientable) such manifolds are the mapping tori of diffeomorphisms
commuting with an Anosov one, namely, mapping tori of Anosov
diffeomorphisms, $\mathbb{T}^3$, and the mapping torus of $-id$
where $id$ is the identity map of $\mathbb{T}^2$.\par
If this conjecture is true, then the fact of being partially
hyperbolic will automatically imply ergodicity in many manifolds. We
prove this for a family of manifolds:

\begin{theorem}\label{nilergodic}
Let $f:N\rightarrow N$ be a conservative partially hyperbolic $C^2$
diffeomorphism where $N\neq \mathbb{T}^3$ is a compact 3-dimensional
nilmanifold. Then, $f$ is ergodic.
\end{theorem}

Recall that a 3-dimensional nilmanifold is a quotient of the
Heisenberg group  (upper triangular $3\times 3$-matrices with ones
in the diagonal) by a discrete subgroup. Sacksteder \cite{sa} proved
that certain affine diffeomorphisms of nilmanifolds are ergodic.
These examples are partially hyperbolic.

Some of our results apply to other manifolds and we obtain, for
instance, that every conservative partially hyperbolic
diffeomorphism of $\mathbb{S}^3$ (or $\mathbb{S}^1\times
\mathbb{S}^2$) is ergodic but is this probably a theorem about the
empty set.

The structure of the proof goes as follows: first of all the results
in \cite{bw,rru} imply that it is enough to prove, in order to
obtain ergodicity, that $f$ satisfies the accessibility property.
Again in \cite{rru} is proved that, in our setting, accessibility
classes are either open or part of a lamination. In Section
\ref{section.lamination} we show that either there are periodic tori
with Anosov dynamics or this lamination extends to a true foliation
without compact leaves and we show that this last case is impossible
for nilmanifolds. In Section \ref{section.invtori} we prove that
such a torus is incompressible (in fact this result is proved in a
more general setting). Since 3-dimensional nilmanifolds do not
support invariant tori with Anosov dynamics we arrive to a
dichotomy: either $f$ is accessible or the accessibility classes are
leaves of a codimension one foliation of $M$ without compact leaves.
On one hand, results of Plante and Roussarie about codimension one
foliations of three manifolds plus its $f$-invariance imply that the
foliation is given by ``parallel'' cylinders (Section
\ref{section.nilergodic}). On the other hand, we exploit the
exponential growth of (Section \ref{section.growth}) unstable curves
to obtain that $f$ should be semiconjugated to a two dimensional
Anosov diffeomorphism (Section \ref{section.heissen}). These two
last facts lead us to a contradiction (Section
\ref{section.nilergodic}).

 All the results of the paper, but the ones that are specific for
 nilmanifolds, can be summarized in the following theorem. Recall
 that accessibility implies ergodicity (see Section
 \ref{section.prelim}).

 \begin{theorem}
 Let $f:M\rightarrow M$ a partially hyperbolic diffeomorphisms of
 an orientable 3-manifold $M$. Suppose that $E^\sigma$ are also
 orientable, $\sigma=s,c,u$, and that $f$ is not accessible. Then
 one of the following possibilities holds:
 \begin{enumerate}
 \item $f$ has a periodic incompressible torus tangent to
 $E^s\oplus E^u$.
 \item \label{lamin} $f$ has an invariant lamination $\Gamma(f)$, tangent to
 $E^s\oplus E^u$, that
 trivially extends to a foliation without compact leaves of $M$. Moreover, the leaves of
 the accessible boundary are periodic, have Anosov dynamics and
 periodic points are dense.
 \item there is a Reebless foliation tangent to $E^s\oplus E^u$.
 \end{enumerate}
  \end{theorem}

 The assumption on the orientability of the bundles and $M$ is not
 essential, in fact, it can be achieved by a finite covering. We
 do not know any example satisfying (\ref{lamin}) in the theorem above (see  Question \ref{franks}).
\medskip

\emph{Acknowledgments:} We are grateful with S. Matsumoto, A. Candel
and, specially, with A. Verjovsky for kindly answer us many
questions about topology and foliations in 3-dimensional manifolds
during the ICM2006.
\end{section}

%% file: prelim2311.tex
\section{Preliminaries}\label{section.prelim}
Let $M$ be a compact Riemannian manifold. In what follows we shall
consider a partially hyperbolic diffeomorphism $f$, that is, a
diffeomorphism admitting a non trivial $Tf$-invariant splitting of
the tangent bundle $TM = E^s\oplus E^c \oplus E^u$, such that all
unit vectors $v^\s\in E^\s_x$ ($\s= s, c, u$) with $x\in M$ verify:
$$\|T_xfv^s\| < \|T_xfv^c\| < \|T_xfv^u\| $$
for some suitable Riemannian metric. It is also required that
$\|Tf|_{E^s}\| < 1$ and $\|Tf^{-1}|_{E^u}\| < 1$. \par%
 It is a known fact that there are foliations $\W^\s$ tangent to the distributions $E^\s$
for $\s=s,u$ (see for instance \cite{bp}). A set $X$ will be called
$\s$-saturated if it is a union of leaves of $\W^\s$, $\s=s,u$. The
{\de accessibility class} $AC(x)$ of the point $x\in M$  is the
minimal $s$- and $u$-saturated set containing $x$. The
diffeomorphism $f$ has the {\de accessibility property} if the
accessibility class of some $x$ is $M$. We shall see that there are
manifolds whose topology implies the accessibility property for all
partially hyperbolic diffeomorphisms, then
ergodicity will follow from the result below:%
\begin{theorem}[\cite{bw},\cite{rru}]
If $f$ is a volume preserving partially hyperbolic diffeomorphism
with the accessibility property and $\dim E^c=1$, then $f$ is
ergodic.
\end{theorem}
The leaf of $\W^\s$ containing $x$ will be called $W^\s(x)$, for
$\s=s,u$. The connected component containing $x$ of the intersection
of $W^s(x)$ with a small $\eps$-ball centered at $x$ is the {\de
$\eps$-local stable manifold of $x$}, and is denoted by
$W^s_\eps(x)$.

Suppose that $V$ is an open invariant set saturated by $s$- and
$u$-leaves ($su$-saturated) and call $\partial^c V$ the accessible
(by center curves) boundary of $V$. Let us observe that $\partial^c
V$ is also an $su$-saturated set, which means that it is a union of
accessibility classes. None of these classes can be open so, due to
Proposition A.3. of \cite{rru}, it is laminated by codimension-one
manifolds tangent to $E^s\oplus E^u$. Let us mention that the leaves
of this lamination are $C^1$ (see \cite{di}).
\par%
Observe also that the proof of Proposition A.5 of \cite{rru} shows
in fact that periodic points are dense in the accessibility classes
of $\partial^cV$ endowed with its intrinsic topology. In other
words, periodic points are dense in each plaque of the lamina of
$\partial^c V$.

We will call $U(f)=\{x\in M; AC(x)$ is open$\}$ and
$\Gamma(f)=M\setminus U(f)$. Then the accessibility property of $f$
is equivalent to  $\Gamma(f)=\emptyset$.

%% file: lami2311.tex
\begin{section}{The $su$-lamination $\Gamma(f)$}\label{section.lamination}

In this section we will study the case when $\emptyset
\varsubsetneq\Gamma(f)\varsubsetneq M$.

By taking a finite covering we can suppose that the bundles
$E^\sigma$ ($\sigma=s,c,u$) and $M$ are orientable.

Let $\Lambda\subset \Gamma(f)$ be an invariant sublamination. Then
its accessible boundary from the complement consists of periodic
lamina with Anosov dynamics and density of periodic points (see
Section \ref{section.prelim}). This implies, thanks to the local
product structure, that stable and unstable leaves of periodic
points are dense in each laminae (see \cite{fr}) and, as a
consequence, that the restriction of an iterated of $f$ that fixes
such a laminae is transitive (always with the intrinsic metric).

Suppose that $\Gamma(f)$ has a compact leaf $F$.  If $\dim (M)=3$,
it is clear that $F$ must be a torus. Since all the leaves of the
invariant sublamination $\Lambda=\bigcup_{n\in \mathbb{Z}}f^n(F)$
are compact we obtain by the comments of the previous section that
there exists a periodic torus leaf with an Anosov dynamics. We deal
with this case in section \ref{section.invtori}. In \cite{rru2}
(Theorem 4.11, Problem 22 and commentary below) we have announced
that $\Gamma(f)$ has always a torus leaf. Unfortunately our proof
has a gap and the following question, up to our knowledge, remains
open even in the codimension one case.

\begin{ques} \label{franks}
Let $f:M\rightarrow M$ be an Anosov diffeomorphism on a complete
riemannian manifold $M$. Is it true that if $\Omega(f)=M$ then $M$
is compact?
\end{ques}

Then in this section we will assume that $\Gamma(f)$ has no compact
leaves (in fact, for our purposes it would be sufficient to assume
that there exists an $f$-invariant sublamination of $\Gamma(f)$
without compact leaves).

Let $\Lambda$ be a lamination and let $V$ be a component of
$M\setminus \Lambda$ and call $V_c$ to the completion of $V$. $V_c$
is a manifold with boundary (we can identify $\partial V_c$ with the
leaves of $\Lambda$ accessible from $V$ and, sometimes, we will
identify $V_c$ with its image under the inclusion
$i_V:V_c\rightarrow M$).

If we take a nice enough covering $\mathcal{U}=\{U_1\dots U_n\}$ by
charts of the lamination $\Lambda$ we have three possible types of
components of $i_V^{-1}(U_j)$:
\begin{itemize}
\item $i_V^{-1}(U_j)$ is diffeomorphic to the open three ball, in
other words $U_j\subset V$.
\item $i_V^{-1}(U_j)$ is diffeomorphic to $[0,1)\times D^2$, where
$D^2$ is the open 2-dimensional disk.
\item $i_V^{-1}(U_j)$ is diffeomorphic to $[0,1]\times D^2$.
\end{itemize}

Since $\mathcal{U}$ is finite there are only a finite number of
components of the first two types. The union of all components of
the third type forms a finite number of connected $I$-bundles
($I=[0,1]$). The union of this $I$-bundles that are not relatively
compact in $V_c$ will be called the interstitial region of $V$. The
union of the other components  (a connected and relatively compact
region) will be called the gut of $V$. Of course this definitions
depends on $\mathcal{U}$. For the basic concepts about laminations
see \cite{gaor}.

Observe that, when $\Lambda$ is a sublamination of $\Gamma(f)$,
$E^c$ is naturally defined in $V_c$.

 We will prove that the
complement of $\Lambda$ consists of $I$-bundles. For this end we
need the following proposition.

\begin{proposition}\label{central.int} Let $f$ be a partially hyperbolic diffeomorphism of a 3-manifold
 $M$.  Let $\Lambda\subset \Gamma(f)$ be an nonempty
$f$-invariant sublamination without compact leaves and $V$ a
component of $M\setminus \Lambda$. Then, $E^c$ is uniquely
integrable in $V_c$.
\end{proposition}

\bp Call $\mathcal{I}$ the interstitial region of $V_c$.
 By Remark 3.7 of \cite{rru1} there exists $\delta>0$
satisfying that if the central bundle is not uniquely integrable at
$x$ and $\gamma$  is a center curve through it, there exists
$n_\gamma >0$ such that the length of $f^n(\gamma)$ or
$f^{-n}(\gamma)$ is greater than $\delta$ for all $n>n_\gamma$. Take
the charts of the lamination in such a way that the interstitial
thickness is small enough (at least less than $\delta$). Then, it is
not difficult to see that if $p$ is a periodic point in
$\mathcal{I}$ then $E^c$ is uniquely integrable at $p$ and,
furthermore, at every point of the connected component of
$W^c(p)\cap V_c$ that contains $p$ (observe that the lamina in
$\partial V_c$ are periodic). This is easily obtained because a
maximal center interval inside $V_c$ containing the periodic point
cannot growth greater than the interstitial thickness. In the
non-periodic case the proof go through similar lines. Suppose that
at $x\in \mathcal{I}$ there are two center curves. By projecting
along the strong stable or strong unstable foliations we can suppose
that these curves are contained in a local center stable or local
center unstable manifold. We assume without loose of generality that
both curves are in the local center unstable manifold of $x$. Then,
we have a ``triangle" formed by the two  center curves and a strong
unstable curve. Now take a  small neighborhood $U$ in the interior
of the triangle and very close to one of the vertices different of
$x$. Since the ``normal" direction to $U$ is stable and
$\Omega(f)=M$ there exists a point $y\in U$ that returns very close
to $U$ for an arbitrarily large $n>0$. By observing that $y$ is the
vertex of a triangle that is very near the previous one (and with
$x$ as a vertex) we obtain that if $\gamma$ is the center curve
joining $x$ and $y$, $f^n(\gamma)$, for $n$ large enough,  have
length greater than $\delta$ which is impossible because $f^n(y)\in
\mathcal{I}$.

Now take $z\in \mathcal{I}$. Then the connected component of
$W^c(z)\cap V_c$ containing $z$ (it is unique) has two extremes
$z_i\in F_i$, with $F_i$ leaves of $\Lambda$, $i=1,2$. By taking an
iterate suppose that $F_1$ (and then $F_2$) is invariant by $f$.
Since $f|_{F_1}$ is transitive we have that for an open and dense
set of points of $F_1$ there is a unique center curve joining this
point and a point in $F_2$ (observe that the same is true for
$F_2$).

We will prove that $\mathcal{V}=\cup_{n\in
\mathbb{Z}}f^n(\mathcal{I})$ is dense in $V_c$. If this is true the
argument of Proposition 1.6 of \cite{bw} can be straightforward
adapted to this case implying unique integrability of the center
unstable and center stable bundles and the proposition.

Observe that $\overline{\mathcal{V}}$ is $su$-saturated. $W^s(x),
W^u(x)\subset \mathcal{V}$ for all $x \in W^c(y)$ with $y\in F_1$
with dense backward and forward orbit. The density of these kind of
points implies the $su$-saturation of $\overline{\mathcal{V}}$.

Then, if $\overline{\mathcal{V}}\neq V_c$, $\overline{\mathcal{V}}$
has nonempty boundary. This boundary consists of  (complete)
$su$-leaves. Close to it there are points of $\mathcal{V}$ that have
center curves that go from $F_1$ to $F_2$. Since these center curves
are contained in $\mathcal{V}$ we arrive to a contradiction.
 \ep

 \begin{theorem}\label{ibundle}
 Let $f$ be a partially hyperbolic diffeomorphism of an orientable 3-manifold
 $M$. Suppose that $\emptyset\varsubsetneq \Lambda\subset \Gamma(f)$
 is an orientable and transversely orientable $f$-invariant sublamination
 without compact leaves such that $\Lambda\neq M$. Then, for all component $V$
 of $M\setminus \Lambda$, $V_c$ is an $I$-bundle.
 \end{theorem}

 The proof of Theorem \ref{ibundle} involves  well known techniques of
 codimension one foliations of three manifolds whose presentation exceeds the purpose of this paper.
 The reader may found them, for instance, in \cite{caco} and \cite{hechir}.

 \medskip

 \bp  Proposition \ref{central.int} implies the existence of center
 stable and center unstable foliations $\mathcal{W}^{cs}, \mathcal{W}^{cu}$ in $V_c$. These foliations
 cannot have Reeb components (see, for instance, \cite{bbi}).
 Observe that, despite the non-compactness of $V_c$, the existence
 of a vanishing cycle implies the existence of a Reeb component.
 This is a consequence of the fact that center stable or center
 unstable disks (whose boundary does not intersect $\partial V_c$) cannot
 intersect $\partial V_c$ because there are no compact stable
 manifolds.  Then the exploding disks of the vanishing cycle cannot
 go inside $\mathcal{I}$, which implies that are all contained in
 the gut. This gives us a Reeb component.

 Now use the notation of Proposition \ref{central.int}.
 Since $F_i\setminus \mathcal{I}$ is compact and
stable and unstable manifolds of periodic points of $f$ in $F_i$ are
dense, $K_i=F_i\setminus \cup_{n\in \mathbb{Z}}f^n(\mathcal{I})$ is
a compact totally disconnected set, $i=1,2$. Take a sequence of
closed curves in $F_1$, $C_n\subset \mathcal{V}$, converging to a
point $p\in K_1$. By taking the other extremes of the center curves
we obtain a sequence of curves in $F_2$ that are also in the gut (it
is clear that they cannot be in $\mathcal{I}$). Taking a convergent
subsequence we have that the limit must be a point of $K_2$
otherwise $K_2$ would contain a nontrivial connected set. Then, for
$n$ large enough, its image in $F_2$, $D_n$, bounds a disk. Consider
the sphere formed by the center curves beginning in $C_n$ and ending
in $D_n$ and the two diks of $F_1$ and $F_2$ bounded by these two
curves. This sphere is homotopically trivial because if this not
were the case we would obtain a vanishing cycle. Then it bounds a
region with trivial fundamental group (in fact, it bounds a ball).
This implies that a center curve starting in $K_1$ must cut $K_2$
(if not we obtain again a vanishing cycle). Then $V_c$ is an
$I$-bundle.
 \ep

 \begin{obs} \label{compact.leaf}
 Observe that Theorem \ref{ibundle} implies that we can complete
 $\Lambda$ to a foliation of $M$ without compact leaves (in particular, Reebless)
 This implies that the inclusion injects the fundamental group of any
 leaf in the fundamental group of $M$. If $M$ is a nilmanifold then
 the fundamental group (of any and then, in particular) of a boundary leaf (periodic with Anosov dynamics and with density of periodic
 points) is isomorphic to a subgroup of a nilpotent group; then itself is nilpotent.
  Since the leaf is an orientable surface this implies that
 it is either a sphere, or a cylinder, or a torus. It is not difficult to prove that
 niether a sphere nor a
 cylinder  support such a dynamics. Then, if $M$ is a
 3-dimensional
 nilmanifold and $\emptyset\varsubsetneq \Gamma(f)\varsubsetneq M$
 then $\Gamma(f)$ has a periodic torus leaf with Anosov dynamics.
 In next sections we will see that this situation is also
 impossible unless the nilmanifold $M$ is $\mathbb{T}^3$.
 \end{obs}

\end{section}

%% file: invtori2311.tex
\begin{section}{Invariant tori}\label{section.invtori}

In this section we show that invariant tori for a diffeomorphism $f$
 with restricted dynamics homotopic to Anosov are incompressible. Recall that a
 two-sided embedded closed surface $S\subset M^3$ other than the sphere is
 {\de incompressible} if and only if
the homomorphism induced by the inclusion
$i_\#:\pi_1(S)\hookrightarrow\pi_1(M)$ is injective; or,
equivalently, by the Loop Theorem, if there is no embedded disc
$D^2\subset M$ such that $D\cap S=\partial D$ and $\partial D\nsim
0$ in $S$ (see, for instance, \cite{hat}).

\begin{theorem}\label{invariant.tori}
Let $f$ be a diffeomorphism of a three-dimensional orientable closed
manifold $M$ and $T\subset M$ an $f$-invariant embedded torus. If
$f|_T$ is isotopic to Anosov, then $T$ is incompressible.
\end{theorem}
\begin{proof}
Let $D^2\subset M$ be an embedded disc such that $D\cap T=\partial
D$ and $\partial D\nsim 0$ in $T$.
 Then, by splitting
$M$ along $T$ we obtain a manifold with boundary $\overline M$ such
that $\partial\overline M=T_1\cup T_2$ where $T_i$, $i=1,2$, are two
tori and at least one of them, say $T_1$, verifies that the
homomorphism induced by the inclusion
$i_\#:\pi_1(T_1)\hookrightarrow\pi_1(\overline M)$ is not
injective.\par
 The diffeomorphism $f$ naturally induces a diffeomorphism $\overline f:\overline M\rightarrow \overline
 M$ fixing  $T_1$ (take $f^2$ if necessary) and such
 that $(\overline f|_{T_1}) _\#:\pi_1(T_1)\rightarrow \pi_1(T_1)$ is a hyperbolic linear automorphism. Moreover, $(\overline f|_{T_1}) _\#$
 leaves $\ker(i_\#)$ invariant.\par
Now, since the eigenspaces of $(\overline f|_{T_1})
 _\#$ have irrational slope, it is not difficult to find
 $(j,0),(0,k)\in \ker(i_\#)$ such that $j,k\in\N\setminus\{0\}$. Let $\alpha, \beta$ be two simple
 closed curves in $T_1$ such that $\alpha$ and $\beta$ meet only
 at a single point and $\alpha^j$ is in the class $(j,0)$ and
 $\beta^k$ in the class $(0,k)$. If we delete all $\partial\overline
 M$ but a small tubular neighborhood of $\alpha$, the Loop Theorem
 gives us a disc $D_\alpha$ embedded in $\overline M$ with boundary
 $\alpha$ and in the same way we obtain $D_\beta$ with boundary
 $\beta$ (this implies $(1,0),(0,1)\in \ker(i_\#)$ and, obviously, $\ker(i_\#)=\pi_1(T_1)$)
 Since we can assume that $D_\alpha$ and $D_\beta$ are transversal,
 this leads to a contradiction with the fact that the intersection
 between $\alpha$ and $\beta$ consist of one point.

\end{proof}
%
%

\end{section}

%% file: growth2311.tex
\begin{section}{Growth of curves}\label{section.growth}

Let us recall the following consequence of Novikov's and Reeb
stability theorems (see \cite{rous}).

\begin{theorem}\label{novikov}
Let $\F$ be a codimension one foliation of a compact 3-manifold $M$.
 and let $\tilde\F$ be the lift of $\F$ to the universal covering
of $M$, $\tilde M$. If $\F$ is Reebless, then either $M$ is
$S^2\times S^1$ and $\F$ is the product foliation or $\tilde M=\R^3$
and $\tilde\F$ is a foliation by planes. Moreover, in this last
case, there is $\epsilon_0>0$ such that $L\cap B_{\epsilon_0}(x)$ is
connected for every $x\in\tilde M$, for every leaf $L$ of
$\tilde\F$.
\end{theorem}
We shall use Theorem \ref{novikov} following the idea in \cite{bbi},
but in our setting. Given  a compact manifold $M$ and $x\in\tilde
M$, let us define $\vol_x(r)=vol\left(B(x,r)\right)$. Notice that
there is $C>0$ such that $\vol_x(r)\leq C\vol_y(r)$ for any two
point $x$ and $y$. So let us fix $x_0\in\tilde M$ and call
$\vol(r)=\vol_{x_0}(r)$.
\begin{proposition}
Let $f:M\to M$ be a partially hyperbolic diffeomorphism of a three
dimensional manifold. Assume that either $E^s\oplus E^u$ or
$E^c\oplus E^u$ is tangent to an invariant foliation $\F$. Then
there is  a constant $C>0$ such that if $I\subset \tilde M$ is an
unstable arc then $\length(I)\leq C \vol\left(\diam(I)\right)$.
\end{proposition}

\bp First of all, observe that $\F$ is a Reebless foliation. This is
a consequence of Proposition 2.1 of \cite{bbi} (see also Theorem
\ref{invariant.tori}) because the boundary of a Reeb component
should be periodic.

 Given an $\varepsilon>0$ small enough
($\varepsilon<\varepsilon_0$ given by Theorem \ref{novikov} and such
that the disk of radius $\varepsilon$ in a leaf is contained in a
chart of the unstable foliation) there exists $\delta>0$ verifying
that the distance between the extremes of an unstable arc of length
greater than $\delta$ is greater than $\varepsilon$. If this were
not the case we have two possibilities:

\begin{enumerate}
\item The two extremes are in different plaques of $\tilde\F$: this
would imply the existence of a closed transversal in $\tilde M$
contradicting that $\F$ is Reebless.

\item The two extremes are in the same plaque of $\tilde\F$: this would
imply that the leaf containing $I$ has nontrivial homotopy
contradicting again that $\F$ is Reebless (by closing $I$ with a
stable or a center arc we obtain an essential closed curve in the
leaf).
\end{enumerate}

Given an unstable arc $I$ there exist at least $\length(I)/2\delta$
disjoint 3-balls of radius $\varepsilon /2$ with center in a point
of $I$. Clearly, the union of these balls is contained in
$B(x,\diam(I))$ for some $x$ which implies:

$$\length(I) \leq \frac{2\delta}{\min \{\vol_y(\varepsilon/2);y\in \tilde M\}}\;\vol_x(\diam(I))\leq C\vol (\diam(I)).$$

\ep

 As nilmanifolds have polynomial growth of volume (of degree at most 4) we get:
\begin{corollary}\label{expgrowth}
In the setting of the proposition above, if $M$ is a nilmanifold
then there is $C>0$ such that $\length(I)\leq C
\left(\diam(I)\right)^4$.
\end{corollary}
\end{section}

%% file: heisen2311.tex
\begin{section}{Nilmanifolds in dimension 3}\label{section.heissen}
Let $\He$ be the group of upper triangular matrices with ones in the
diagonal. This is the non-abelian nilpotent simply connected three
dimensional Lie group. We may identify $\He$ with the pairs $(x,t)$
where $x=(x_1,x_2)\in\R^2$, $t\in\R$,
$(x,t)\cdot(y,s)=(x+y,t+s+x_1y_2)$ and $(x,t)^{-1}=(-x,x_1x_2-t)$.
We have the projection $p:\He\to\R^2$, $p(x,t)=x$ which is also an
homomorphism.

If we denote with $\he$ the Lie algebra of $\He$, then we may also
identify $\he$ with the pairs $(x,t)$ where $x=(x_1,x_2)\in\R^2$,
$t\in\R$. We have the exponential map $\exp:\he\to\He$ given by
$\exp(x,t)=(x,t+\frac{1}{2}x_1x_2)$, $\exp$ is one to one and onto;
and its inverse, the logarithm, $\log:\He\to\he$ is given by
$\log(x,t)=(x,t-\frac{1}{2}x_1x_2)$.

The homomorphisms from $\He$ to $\He$ are of the form
$L(x,t)=\left(Ax,l(x,t)\right)$, where
$$
A=\left(\begin{array}{cc}
a & b\\
c & d
\end{array}\right)\;\;\;\;\mbox{and}\;\;\;\; l(x,y)=\alpha x_1+\beta x_2+\det(A)t+\frac{ac}{2}x_1^2+\frac{bd}{2}x_2^2+bcx_1x_2.
$$

If we denote with $\hat L:\he\to\he$, $\hat L=D_0 L$, $\hat L$ is induced by the matrix
\begin{displaymath}
\hat L=\left(
\begin{array}{ccc}
a & b & 0\\
c & d & 0\\
\alpha &\beta & \det(A)
\end{array}\right)
\end{displaymath}
and $\exp\left(\hat L(x,t)\right)=L\left(\exp(x,t)\right)$.

The centralizer of $\He$ is exactly $\He_1=[\He,\He]$ which consist
the elements of the form $(0,t)$. Any homomorphism from $\He$ to
$\He$ must leave $\He_1$ invariant.

Any lattice in $\He$ is isomorphic to $\Gamma_k=\{(x,y):x\in\Z^2,y\in \frac{1}{k}\Z\}$, for $k$ a positive integer.

The automorphisms of $\He$ are exactly the ones with $\det(A)\neq 0$ and the automorphisms leaving $\Gamma_k$ invariant are the ones with $A\in GL(2,\Z)$
(the matrices with integral entries and determinant $\pm 1$) and $\alpha,\beta\in\frac{1}{k}\Z$. On the other hand, every automorphisms of $\Gamma_k$
extends to an automorphisms of $\He$.
\begin{lema}\label{centralizer}
If $S$ is a subgroup of $\Gamma_k$ isomorphic to $\Z^2$, then $S\cap
\He_1\neq\{(0,0)\}$.
\end{lema}
\begin{proof}
Let $(x,t)$ and $(y,s)$ generate $S$. Then,
$(x,t)\cdot(y,s)=(y,s)\cdot(x,t)$ implies $x_1y_2=x_2y_1$. So,
$y=\frac{p}{q}x$. Now, $(x,t)^p\cdot(y,s)^{-q}=(0,u)$ for some $u\in
\R$. The fact that $(x,t)$ and $(y,s)$ generate $S$ implies $u\ne0$,
so $(0,u)\in\He_1\cap S$
\end{proof} We define the quotient compact nilmanifold $N_k=\He/\Gamma_k$ by
the relation $(x,t)\sim (y,s)$ iff
$(x,t)^{-1}\cdot(y,s)\in\Gamma_k$. The first homotopy group is
$\pi_1(N_k)=\Gamma_k$ and two maps of $N_k$ to itself are homotopic
if and only if their action on $\Gamma_k$ coincide. Moreover, any
map from $N_k$ to itself is homotopic to an automorphism as the
described above leaving $\Gamma_k$ invariant. Given $f:N_k\to N_k$,
with induced automorphism on $\Gamma_k$, $f_\#=L$, if $F:\He\to\He$
is a lift of $f$ to $\He$, then $F(zn)=F(z)L(n)$ for every
$n\in\Gamma_k$, $z\in\He$. Moreover, $F(z)=u(z)L(z)$, where
$u:\He\to\He$ is such that $u(zn)=u(z)$ for every $n\in\Gamma_k$,
$z\in\He$ and for $k>0$, $F^k(z)=u_k(z)L^k(z)$ where
\begin{eqnarray}\label{fk}
u_{k+1}(z)=u\left(F^k(z)\right)L\left(u\left(F^{k-1}(z)\right)\right)\dots
L^k\left(u(z)\right)L^{k+1}(z)
\end{eqnarray}

The projection $p:N_k\to\T^2$,
$p\left((x,t)\cdot\Gamma_k\right)=x+\Z^2$ induces an isomorphism
$p_*:H_1(N_k)\to H_1(\T^2)$. Let $d$ be a right invariant metric on
$\He$, for instance
$d\left((x,t),0\right)=|x|+|t-\frac{1}{2}x_1x_2|$ where
$x=(x_1,x_2)$. The exponential $\exp:\he\to\He$ become an isometry
if the metric in $\he$ is $|(x,t)|=|x|+|t|$.

We have the following proposition:
\begin{proposition}\label{acces.or.int}
Let $f:N_k\to N_k$ be a partially hyperbolic diffeomorphism. If $f$
has not the accessibility property then $E^s\oplus E^u$ integrates
to a foliation such that any closed nonempty $f$-invariant set
saturated by leaves is the whole manifold $M$. Moreover, the action
of $f_*=A$ on $H_1(N_k,\Z)$ is hyperbolic and hence there is a
semiconjugacy $h:N_k\to\T^2$ homotopic to $p$ such that $h\circ
f=Ah$.
\end{proposition}
\begin{proof}
Using Remark \ref{compact.leaf} and Theorem \ref{invariant.tori}, we
have that either $E^s\oplus E^u$ integrates to a foliation with the
minimality property above or there is a torus $T$ invariant by $f^k$
whose homotopy group injects on $\pi_1(N_k)$ and such that
$f^k_\#|\pi_1(T)$ is hyperbolic.

If this last were the case, as $\pi_1(T)\sim\Z^2$, using Lemma
\ref{centralizer} we get that $\pi_1(T)\cap\He_1\neq\{0\}$. On the
other hand, for any automorphism $L$ of $\Gamma_k$,
$L(x,t)=(x,t)^{\pm 1}$ for every $(x,t)\in\He_1$ which gives a
contradiction.\par
So we have that $E^s\oplus E^u$ integrates to a foliation. Let us
call $L=f_\#$. Notice that $L(x,t)=(Ax,l(x,t))$ where $A=f_*$ the
action of $f$ on $H_1(N_k,\Z)$. Thus, if $A$ is not hyperbolic, then
it is not hard to see that for $k>0$, $d(L^kz,0)\leq Ck^2d(z,0)$ for
some constant $C>0$. Take $F(z)=u(z)Lz$ a lift of $f$ to $\He$.
Then, using Formula \ref{fk} and that $d$ is a right invariant
metric, we have that for every $z\in\He$,
\begin{eqnarray*}
d\left(F^k(z),0\right)&\leq&\sum_{i=0}^{k-1}d\left(L^iu\left(f^{k-1-i}(z)\right),0\right)+d\left(L^k(z),0\right)\\
&\leq&d\left(u\left(f^{k-1}(z)\right),0\right)+\sum_{i=1}^{k-1}Ci^2d\left(u\left(f^{k-1-i}(z)\right),0\right)+Ck^2d\left(z,0\right)\\
&\leq&Ck^3+Ck^2d(z,0)
\end{eqnarray*}
since $d(u(z),0)\leq C$ for every $z\in\He$.

Thus, if $I$ is an unstable arc in $\He$, then for $k>0$,
$diam\left(F^k(I)\right)\leq C(I)k^3$ for some constant $C(I)$ that
does not depend on $k$. On the other hand, by Corollary
\ref{expgrowth}, the growth of $F^k(I)$ should be exponential. Thus
we get a contradiction and hence $A$ must be hyperbolic. The
existence of $h$ follows from standard arguments (see \cite{fr}).
\end{proof}

\end{section}

%% file: nilergodic2311.tex
\begin{section}{Proof of Theorem
\ref{nilergodic}}\label{section.nilergodic}

As was shown in Proposition \ref{acces.or.int} if a conservative
partially hyperbolic diffeomorphism $f:N_k\to N_k$ has not the
accessibility property then $E^s\oplus E^u$ integrate to a foliation
$\F^{su}$ with the following minimality property: any closed,
nonempty, $f$-invariant and saturated by leaves set is the whole
manifold $M$. In this last section we shall prove the existence of
such a foliation leads us to a contradiction. Without lose of
generality we may assume, by taking a double covering if necessary,
that $\F^{su}$ is transversely orientable. Observe that the double
covering of a nilmanifold is again a nilmanifold.

Observe that $\F^{su}$ has  no compact leaves. On one hand there are
not periodic compact leaves by the minimality property of the
foliation. On the other hand, if there is a compact noninvariant
leaf $T$ (it must be a torus) then $\overline{\{f^n(T);n\in
\mathbb{Z}\}}=M$. This implies that all the leaves are tori and $M$
is the mapping torus of a linear automorphism of $\mathbb{T}^2$ that
commutes with a hyperbolic one. The only three dimensional
nilmanifold satisfying this is $\mathbb{T}^3$.

We need the following theorem whose proof is essentially in
\cite{plante.tmeas} (see also \cite{hechir}).

\begin{theorem}\label{minfinitos}
Let $\F$ be a codimension one $C^0$-foliation without compact leaves
of a three dimensional compact manifold $M$. If $\pi_1(M)$ has
non-exponential growth the number $k$ of minimal sets that support
an $\F$-invariant measure satisfies $0<k<+\infty$.
\end{theorem}

\bp Theorems 7.3, 4.1 and 6.3 of \cite{plante.tmeas} implies the
existence of a minimal set supporting an $\F$-invariant measure and
Proposition 8.5 of the same paper implies finiteness. Observe that
the $C^1$ hypothesis of Theorem 7.3 comes from Novikov's Theorem
that is valid in the $C^0$ setting (\cite{solodov}). \ep

\begin{lema}
$\F^{su}$ is minimal.
\end{lema}

\bp $\F^{su}$ has not compact leaves and the fundamental group of a
nilmanifold has polynomial growth. Call $K_1,\dots,K_k$ to the
minimal sets given by Theorem \ref{minfinitos}. $K=K_1\cup\dots\cup
K_k$ is an $\F^{su}$-saturated set and it is $f$-invariant. Thus
$K=M$ proving the minimality of $\F^{su}$.\ep

As a consequence of the existence of a invariant transversal
measure, we obtain that the holonomy of $\F^{su}$ is trivial
(Theorem 2.3.1 of Chapter X of \cite{hechir}).

 Since $\F^{su}$ has no compact leaves, an incompressible torus can
 be isotoped in such a way
it becomes everywhere transversal to it (see \cite{rous2}).
Roussarie result is stated with the assumption of some extra
differentiability for the foliation but the general position
arguments on which his result is based are valid for $C^0$
foliations (see \cite{hechir} and, also, \cite{gabai}). Now cutting
along the transverse torus and arguing as in Theorem 4.1 of
\cite{plante.solv} we obtain that the foliation on
$\mathbb{T}^2\times [0,1]$ may be isotoped to a foliation obtained
from a foliation of the torus $\mathcal{P}$ covered by parallel
lines of the plane taking as leaves the surfaces (leaf in
$\mathbb{T}^2$)$\times [o,1]$. In order to obtain the original
nilmanifold we have to identify the two boundary tori by a
homeomorphism isotopic to $\left(\begin{array}{cc}
1 & k\\
0 & 1
\end{array}\right)$. This implies that the lines of $\mathcal{P}$
are parallel to the vector $(1,0)$, thus $\mathcal{P}$ is a
foliation by circles. Moreover, this curves correspond to the
homotopy class of the center $\He_1$.

\begin{obs}\label{ciriso}
Observe that the leaves of $\F^{su}$ are subfoliated by circles that
are (uniformly) isotopic to the corresponding to $x=constant$ where
$(x,y)\in \He$.
\end{obs}

Let $h:N_k\rightarrow \mathbb{T}^2$ be the semiconjugacy given by
Proposition \ref{acces.or.int}.

\begin{lema}
There is $w\in N_k$ such that $h(W^\sigma (w))=W^\sigma (h(w))$ with
$\sigma =u$ or $s$.
\end{lema}

\bp First of all we claim that there exists $x$ such that
$h(W^u(x))$ or $h(W^u(x))$ contains more than one point. If for all
$y\in N_k$ we have that $h(W^u(y))=h(W^s(y))=q_y\in \mathbb{T}^2$
then $h(AC(x))=q_x$ and the minimality of $\F^{su}$ implies that
$AC(x)$ is dense and $h(N_k)=q_x$ contradicting the sobrejectivity
of $h$.

Take $x\in N_k$ such that $h(W^u(x))$ is a nontrivial interval of
$W^u(h(x))$ and $z\in W^u(x)$ such that $h(z)$ is an interior point
of $h(W^u(x))$. A standard argument shows that any point $w\in
\omega(z)$ satisfies that $h(W^u(w))=W^u(h(w))$. \ep

Take a point $w$ satisfying the lemma (suppose that it works in the
$u$ case) and call $F$ to the $\F^{su}$-leaf through $w$. There
exists an injective immersion $i:\mathbb{R}\times
\mathbb{S}^1\rightarrow N_k$ such that $i(\mathbb{R}\times
\mathbb{S}^1)=F$ and consider $j:h\circ i:\mathbb{R}\times
\mathbb{S}^1\rightarrow \mathbb{T}^2$. The previous considerations
about $\F^{su}$ implies that $j_\#=h_\#\circ i_\#=p_\#\circ
i_\#:\pi_1(\mathbb{R}\times \mathbb{S}^1)\rightarrow
\pi_1(\mathbb{T}^2)$ is trivial. Then there exists
$\tilde{\j}:\mathbb{R}\times \mathbb{S}^1\rightarrow \mathbb{R}^2$
such that $j=\pi\circ \tilde{\j}$ where $\pi:\mathbb{R}^2\rightarrow
\mathbb{T}^2$ is a covering projection.

Remark \ref{ciriso} implies that $\tilde{\j}(\mathbb{R}\times
\mathbb{S}^1)$ is contained in bounded neighborhood of
$W^u(\tilde{\j}(i^{-1}(w))\subset \tilde{\j}(\mathbb{R}\times
\mathbb{S}^1)$.

Consider now $i^{-1}(W^s(w))$. Since the foliation $\F^{s}$ has
neither singularities nor compact leaves, $i^{-1}(W^s(w))$ is
unbounded. It is not difficult to see that this implies that
$\tilde{\j}(i^{-1}(W^s(w)))$ is unbounded but it is at bounded
distance of $W^u(\tilde{\j}(i^{-1}(w)))$ which is a contradiction.
This proves Theorem \ref{nilergodic}.
\end{section}